# Об эффективности одного метода рандомизации зеркального спуска в задачах онлайн оптимизации[1]

*А.В. Гасников, Ю.Е. Нестеров, В.Г. Спокойный*

gasnikov@yandex.ru, yurii.nesterov@uclouvain.be, spokoiny@wias-berlin.de

Лаборатория структурных методов анализа данных в предсказательном моделировании,

Кафедра математических основ управления,

Факультет управления и прикладной математики МФТИ;

Институт Проблем Передачи Информации РАН;

НИУ Высшая школа экономики

В статье предложена рандомизированная онлайн версия метода зеркального спуска. Отличие от имеющихся версий заключается в способе рандомизации. Рандомизация вводится не на этапе вычисления субградиента функции, как это повсеместно принято, а на этапе проектирования на единичный симплекс. В результате получается покомпонентный субградиентный спуск со случайным выбором компоненты, допускающий онлайн интерпретацию. Это наблюдение, например, позволило единообразно проинтерпретировать результаты о взвешивании экспертных решений, предложить наиболее эффективный способ поиска равновесия в антагонистической матричной игре с разреженной матрицей.

**Ключевые слова:** метод зеркального спуска, метод двойственных усреднений, онлайн оптимизация, экспоненциальное взвешивание, многорукие бандиты, взвешивание экспертов, стохастическая оптимизация, рандомизация.

---

# 1. Введение

В конце 70-х годов А.С. Немировским и Д.Б. Юдиным был предложен итерационный метод решения негладких выпуклых задач оптимизации (см., например, [1]), который можно интерпретировать как разновидность метода проекции субградиента, когда проектирование понимается, например, в смысле расстояния Брэгмана (Кульбака–Лейблера) [2], или как прямодвойственный метод [3]. Этот метод, получивший название *метода зеркального спуска* (МЗС), позволяет хорошо учитывать структуру множества, на котором происходит оптимизация (например, симплекса). Как и многие другие методы решения негладких выпуклых оптимизационных задач, этот метод требует $O(M^2R^2/\varepsilon^2)$ итераций, где $\varepsilon$ – точность найденного решения по функции, что соответствует нижним оценкам по $\varepsilon$ [1]. Однако константа $M$, равномерно ограничивающая норму субградиента оптимизируемой функции, и размер множества $R$ зависят от выбора нормы в пространстве, в котором ведется оптимизация. Так, если мы выбрали норму в нашем пространстве $l_p$ ($1 \le p \le \infty$), то $M$ – есть сопряженная $l_q$-норма субградиента ($1/p + 1/q = 1$), а $R^2$ – есть "размер" множества в "метрике" сильно выпуклой относительно $l_p$, с константой сильной выпуклости $\alpha \ge 1$.

**Замечание 1.** Слово "размер" взято в кавычки, потому что в действительности то, что задается, мы интерпретируем в данном контексте как квадрат размера (физически правильнее квадратом размера называть $R^2/\alpha$), поскольку "метрика" сильно выпуклая относительно нашей "рабочей" нормы в этом пространстве. Слово "метрика" взято в кавычки, потому что может быть не выполнено одно свойств метрики – нет симметричности, например, для расстояния (дивергенции) Брэгмана.

При таких предположениях говорят, что выбранная "метрика" порождает прокс-структуру на множестве. Например, когда множество, на котором происходит оптимизация является симплексом, то, как правило, выбирают $p = 1$, а "метрику" задают расстоянием Брэгмана. При этом "проекция" на симплекс согласно такому расстоянию считается по явным формулам (экспоненциальное взвешивание). В работе [4] была выдвинута гипотеза о том, что применительно к задачам стохастической оптимизации на единичном симплексе (не онлайн) такой выбор нормы и расстояния являются наилучшими с точки зрения зависимости $M^2R^2$ от размера пространства $n$ (в типичных приложениях эта зависимость $\sim \ln n$). Однако в определенных ситуациях (в задачах о многоруких бандитах, когда $M^2R^2 \sim n\ln n$) удается выиграть логарифмический по $n$ фактор (см. пример 1 п. 4), более подходящим образом выбирая расстояние [5]. При этом теряется возможность явного вычисления проекции.

В работах [3] – [8] исследовались стохастические (рандомизированные) версии МЗС. В том числе и онлайн [5]. При этом анализировалась ситуация, когда именно градиент функции выдается оракулом со случайными шумами, но несмещенным образом. Такая релаксация детерминированного МЗС оказалась весьма полезной применительно к задачам адаптивного агрегирования оценок [4], оптимизации в пространствах огромных размеров [7], [8], задачах о многоруких бандитах и т.п. [5], [9] – [11].

В работе [1] была также отмечена возможность онлайн интерпретации МЗС. Впоследствии, у разных авторов можно найти заметки на эту тему [3], [5], [9], [11] – [16]. Наблюдение состоит в том, что ничего не изменится с точки зрения изучения сходимости метода (и его стохастической версии), если на каждом шаге допускать, что функция меняется, причем, возможно, враждебным образом (при этом оставаясь в классе выпуклых функций с ограниченной нормой субградиента).

В данной работе приводятся две версии стохастического онлайн МЗС. Первая версия, более менее, классическая. Приблизительно в таком же виде её уже можно было встретить в литературе у разных авторов [3] – [5], [7] – [16]. Точнее говоря, предложенная версия аккумулирует в себе в виде частных случаев многие известные ранее версии МЗС. Вторая неявно была предложена в работе [6] применительно к поиску равновесия в антагонистической матричной игре (онлайн модификация в [6] не была затронута, равно как и связь предложенного метода с МЗС). Согласно работе [6], мы рандомизируем не на этапе вычисления субградиента функции, как это общепринято [7], [8], а на этапе проектирования на симплекс. В результате получается покомпонентный субградиентный спуск со случайным выбором компоненты, который, как будет ниже показано, допускает онлайн интерпретацию. Получив такой метод, мы расширяем множество тех задач онлайн оптимизации, к которым можно применять МЗС.

## 2. Онлайн МЗС со стохастическим субградиентом

Рассмотрим задачу стохастической онлайн оптимизации (запись $E_{\xi^k}\left[f_k(x;\xi^k)\right]$ означает, что математическое ожидание берется по $\xi^k$, то есть $x$ и $f_k$ понимаются в такой записи не случайными)

$$\frac{1}{N}\sum_{k=1}^{N}E_{\xi^k}\left[f_k(x;\xi^k)\right] \to \min_{x \in S_n(1)}, \ S_n(1) = \left\{x \geq 0: \ \sum_{i=1}^{n}x_i = 1\right\}, \quad (1)$$

при следующих **условиях**:

1. $E_{\xi^k}\left[f_k(x;\xi^k)\right]$ – выпуклые функции (по $x$), для этого достаточно выпуклости по $x$ функций $f_k(x;\xi^k)$ и независимости распределения $\xi^k$ от $x$;

2. Существует такой вектор $\nabla_x f_k(x;\xi^k)$, который для компактности будем называть субградиентом, хотя последнее верно не всегда (см. пример 1 п. 4), что

$$E_{\xi^k}\left(\nabla_x f_k(x;\xi^k) - \nabla_x E_{\xi^k}\left[f_k(x;\xi^k)\right]\Big|\Xi^{k-1}\right) \equiv 0,$$

где $\Xi^{k-1}$ – $\sigma$-алгебра, порожденная случайными величинами $\xi^1,...,\xi^{k-1}$. Далее везде в статье мы будем использовать обозначения обычного градиента для векторов, которые мы назвали здесь субградиентами. В частности, если мы имеем дело с обычным субградиентом, то запись $\nabla_x f_k(x;\xi^k)$ в вычислительном контексте (например, в итерационной процедуре МЗС, описанной ниже) означает какой-то его элемент (не важно какой именно), а если в контексте проверки условий (например, в условии 3 ниже), то $\nabla_x f_k(x;\xi^k)$ пробегает все элементы субградиента (говорят также, субдифференциала);

3. $\left\|\nabla_x f_k(x;\xi^k)\right\|_\infty \leq M$ – (равномерно, с вероятностью 1) ограниченный субградиент.

Для справедливости части утверждений достаточно требовать одно из следующих (более слабых) условий:

а) $E_{\xi^k}\left[\left\|\nabla_x f_k(x;\xi^k)\right\|_\infty^2\right] \leq M^2$;   б) $E_{\xi^k}\left[\exp\left(\frac{\left\|\nabla_x f_k(x;\xi^k)\right\|_\infty^2}{M^2}\right)\bigg|\Xi^{k-1}\right] \leq \exp(1)$.

Задача (1) является лишь компактной (и далеко не полной) записью настоящей постановки задачи стохастической онлайн оптимизации. В действительности, требуется подобрать последовательность $\{x^k\}$ (в п. 2 $\{x^k\} \in S_n(1)$, а в п. 3 и ряде примеров п. 4 при дополнительном ограничении, что $\{x^k\}$ выбираются с возможными повторениями среди вершин симплекса $S_n(1)$) исходя из доступной исторической информации ($x^k$ может зависеть только от $\{x^1,\xi^1,f_1(\cdot);...;x^{k-1},\xi^{k-1},f_{k-1}(\cdot)\}$) так, чтобы минимизировать псевдо регрет:

$$\frac{1}{N}\sum_{k=1}^N E_{\xi^k}\left[f_k(x^k;\xi^k)\right] - \min_{x \in S_n(1)}\frac{1}{N}\sum_{k=1}^N E_{\xi^k}\left[f_k(x;\xi^k)\right]$$

или регрет

$$E_{\xi^1,...,\xi^N}\left[\min_{x \in S_n(1)}\frac{1}{N}\sum_{k=1}^N\left(f_k(x^k;\xi^k) - f_k(x;\xi^k)\right)\right].$$

В данном пункте мы сосредоточимся на минимизации псевдо регрета на основе информации $\{\nabla f_1(x^1;\xi^1);...;\nabla f_{k-1}(x^{k-1};\xi^{k-1})\}$ при расчете $x^k$. Онлайновость постановки задачи допускает, что на каждом шаге $k$ функция $f_k$ может подбираться из рассматриваемого класса функций враждебно по отношению к используемому нами методу генерации последовательности $\{x^k\}$. В частности, в этом пункте $f_k$ может зависеть от $\{x^1,\xi^1,f_1(\cdot);...;x^{k-1},\xi^{k-1},f_{k-1}(\cdot);x^k\}$.

В целом в данной статье мы ограничиваемся рассмотрением задач минимизации псевдо регрета, когда на каждом шаге мы можем получить независимую реализацию стохастического субградиента в одной указанной нами (допустимой) точке. Приведенные в статье результаты можно распространить на случай, когда градиент выдается не точно (с не случайной ошибкой), выдается не полностью (скажем, вместо градиента выдается производная по выбранному направлению) или вместо градиента выдается только значение функции [14]. Впрочем, немного об этом написано далее (пример 1 п. 4). Другим способом релаксации исходной постановки является возможность несколько раз обращаться на одном шаге за значением градиента функции и(или) значением самой функции и взвешивать $f_k(x;\xi^k)$ в (1) разными весами (см. [14], [17]). Результаты статьи также можно распространить и на случай, когда функции $E_{\xi^k}\left[f_k(x;\xi^k)\right]$ равномерно сильно выпуклые по $x$ с константой $\mu$ [5], [11], [13], [16]. При этом выбирается евклидова прокс-структура, поскольку в сильно выпуклом случае (в стохастическом и не стохастическом) игра на выборе прокс-структуры не может дать выгоды. Неулучшаемая и достижимая оценка в этом случае будет иметь следующий вид $\varepsilon = \mathrm{O}\left(M^2 \ln N / (\mu N)\right)$. В не онлайн стохастическом случае эта оценка (с точностью до фактора $\ln N$) также будет неулучшаемой. Отметим при этом, что (для рассмотренных выпуклых и сильно выпуклых задачах) в отличие от не онлайн случая, в онлайн случае игра на гладкости функций $E_{\xi^k}\left[f_k(x;\xi^k)\right]$ и(или) отсутствии стохастичности ($f_k(x;\xi^k) \equiv f_k(x)$) не дает никаких дивидендов (выписанные нижние оценки сохранятся).

Можно обобщить приведенную постановку и последующие результаты на задачи композитной оптимизации [18]. Если функция, которая добавляется ко всем $f_k$, – линейная, то мы даже можем ее не знать (просто знать, что она есть и одна и та же). Тогда опустив первое слагаемое в сумме (1) и переписав условие 3 в виде

$$\left\|\nabla_x f_k(x;\xi^k) - \nabla_x f_m(y;\xi^m)\right\|_\infty \le M, \text{ для всех } x, y \in S_n(1),\ k, m \in \mathbb{N},$$

можно перенести результаты статьи [19] на такой онлайн контекст.

Выше мы исходили из того, что оптимизация ведется на единичном симплексе. Возникает вопрос: насколько все, что приведено в статье, обобщается на более общий случай? Собственно говоря, ответ на этот вопрос частично известен уже давно [1]. Приведенные в п. 2 рассуждения универсальны, то есть если исходить из оптимизации на каком-нибудь другом выпуклом компакте (от условия компактности (ограниченности) множества можно отказаться [3], поскольку, в действительности, в оценку числа итераций входит не размер множества, а "расстояние" от точки старта до решения), то задав норму в прямом пространстве и расстояние (сильно выпуклое относительно этой нормы), согласно которому будет осуществляться проектирование субградиента на этот компакт, можно повторить аналогичные рассуждения [20].

Для решения задачи (1) воспользуемся адаптивным методом зеркального спуска (точнее двойственных усреднений) в форме [3], [4]. Положим $x_i^1 = 1/n$, $i = 1,...,n$. Пусть $t = 1,..., N-1$.

**<u>Алгоритм МЗС1-адаптивный / Метод двойственных усреднений</u>**

$$x_i^{t+1} = \frac{\exp\left(-\frac{1}{\beta_{t+1}}\sum_{k=1}^{t}\frac{\partial f_k\left(x^k;\xi^k\right)}{\partial x_i}\right)}{\sum_{l=1}^{n}\exp\left(-\frac{1}{\beta_{t+1}}\sum_{k=1}^{t}\frac{\partial f_k\left(x^k;\xi^k\right)}{\partial x_l}\right)}, \quad i = 1,...,n, \quad \beta_t = \frac{M\sqrt{t}}{\sqrt{\ln n}}.$$

Не сложно показать, что этот метод представим также в виде:

$$x^{t+1} = \arg\min_{x \in S_n(1)}\left\{\sum_{k=1}^{t}\left\{f_k\left(x^k;\xi^k\right) + \left\langle\nabla f_k\left(x^k;\xi^k\right), x - x^k\right\rangle\right\} + \beta_{t+1}V(x)\right\}$$

или

$$\begin{cases} y^k = y^{k-1} - \gamma_k \nabla_x f_k\left(x^k;\xi^k\right) \\ x^{k+1} = \nabla W_{\beta_{k+1}}\left(y^k\right) \end{cases}, \ y^0 = 0, \ \gamma_k \equiv 1, \ \beta_k = \frac{M}{\sqrt{\ln n}}\sqrt{k}, \ k = 1,...,N, \quad (2)$$

где

$$W_\beta(y) = \sup_{x \in S_n(1)}\left\{\langle y, x\rangle - \beta V(x)\right\} = \beta\ln\left(\frac{1}{n}\sum_{i=1}^{n}\exp(y_i/\beta)\right),$$

$$V(x) = \ln n + \sum_{i=1}^{n} x_i \ln x_i.$$

Рассуждая подобно [3], [4], [7], можно получить следующий результат.

**Теорема 1.** *Пусть справедливы условия 1, 2, 3.а, тогда*

$$\frac{1}{N}\sum_{k=1}^{N}E\left[f_k\left(x^k;\xi^k\right)\right] - \min_{x \in S_n(1)}\frac{1}{N}\sum_{k=1}^{N}E_{\xi^k}\left[f_k\left(x;\xi^k\right)\right] \leq 2M\sqrt{\frac{\ln n}{N}}.$$

*Если $f_k \equiv f$, а $\{\xi^k\}$ – независимы и одинаково распределены, как $\xi$, то*

$$E\left[f\left(\frac{1}{N}\sum_{k=1}^{N}x^k;\xi\right)\right] - \min_{x \in S_n(1)}E_\xi\left[f(x;\xi)\right] \leq 2M\sqrt{\frac{\ln n}{N}}.$$

*Пусть справедливы условия 1, 2, 3, тогда при $\Omega \geq 0$*

$$P_{x^1,...,x^N}\left\{\frac{1}{N}\sum_{k=1}^{N}E_{\xi^k}\left[f_k\left(x^k;\xi^k\right)\right] - \min_{x \in S_n(1)}\frac{1}{N}\sum_{k=1}^{N}E_{\xi^k}\left[f_k\left(x;\xi^k\right)\right] \geq \frac{2M}{\sqrt{N}}\left(\sqrt{\ln n} + \sqrt{8\Omega}\right)\right\} \leq \exp(-\Omega).$$

*Если $f_k \equiv f$, а $\{\xi^k\}$ – независимы и одинаково распределены, как $\xi$, то*

$$P_{x^1,\ldots,x^N}\left\{ E_\xi\left[ f\left(\frac{1}{N}\sum_{k=1}^{N} x^k;\xi\right)\right] - \min_{x\in S_n(1)} E_\xi\left[ f(x;\xi)\right] \geq \frac{2M}{\sqrt{N}}\left(\sqrt{\ln n} + \sqrt{8\Omega}\right)\right\} \leq \exp(-\Omega).$$

**Замечание 2.** Запись

$$"P_{x^1,\ldots,x^N}\left\{ \frac{1}{N}\sum_{k=1}^{N} E_{\xi^k}\left[ f_k(x^k;\xi^k)\right] - \ldots"$$

означает, что под вероятностью мы считаем математическое ожидание по $\xi^k$, которое, вообще говоря, зависит и от $\xi^1,\ldots,\xi^{k-1}$ (мы не предполагаем независимости $\{\xi^k\}$), как бы "замораживая" (считая не случайными) $x^k$, то есть забывая про то, что $x^k$ тоже зависит от $\xi^1,\ldots,\xi^{k-1}$. А вероятность берется как раз по $\{x^k\}$, с учетом того, что такая зависимость есть (см. определение алгоритма МЗС1).

**Доказательство теоремы 1.** Возьмем за основу обозначения МЗС1, приведенные в (2). Рассуждая далее аналогично [3], [4], [7], получим (заметим, что $\beta_{k+1} \geq \beta_k > 0$):

$$W_{\beta_{k+1}}(y^k) \leq W_{\beta_k}(y^k) = W_{\beta_k}(y^{k-1}) + \int_0^1 \left(y^k - y^{k-1}\right)^T \nabla W_{\beta_k}\left(\tau y^k + (1-\tau) y^{k-1}\right) d\tau =$$

$$= W_{\beta_k}(y^{k-1}) - \gamma_k \nabla_x f_k(x^k;\xi^k)^T \nabla W_{\beta_k}(y^{k-1}) -$$

$$-\gamma_k \nabla_x f_k(x^k;\xi^k)^T \int_0^1 \left(\nabla W_{\beta_k}\left(\tau y^k + (1-\tau) y^{k-1}\right) - \nabla W_{\beta_k}(y^{k-1})\right) d\tau \overset{(*)}{\leq}$$

$$\overset{(*)}{\leq} W_{\beta_k}(y^{k-1}) - \gamma_k \nabla_x f_k(x^k;\xi^k) \nabla W_{\beta_k}(y^{k-1}) +$$

$$+ \gamma_k \left\|\nabla_x f_k(x^k;\xi^k)\right\|_\infty \int_0^1 \left\|\nabla W_{\beta_k}\left(\tau y^k + (1-\tau) y^{k-1}\right) - \nabla W_{\beta_k}(y^{k-1})\right\|_1 d\tau \leq$$

$$\leq W_{\beta_k}(y^{k-1}) - \gamma_k \nabla_x f_k(x^k;\xi^k) \nabla W_{\beta_k}(y^{k-1}) + \frac{\gamma_k^2 \left\|\nabla_x f_k(x^k;\xi^k)\right\|_\infty^2}{2\alpha\beta_k},$$

последнее неравенство следует из того, что [3], [4], [7]:

$$\left\|\nabla W_\beta(\tilde{y}) - \nabla W_\beta(y)\right\|_1 \leq \frac{1}{\alpha\beta}\|\tilde{y} - y\|_\infty,$$

где $\alpha = 1$ – константа сильной выпуклости $V(x)$ в 1-норме. Мы специально выделили (*) неравенство, которое иногда (например, в задачах о многоруких бандитах) бывает довольно грубым. В работах [5], [14] указан способ уточнения этого неравенства.

Суммируя эти неравенства, учитывая, что $\nabla W_{\beta_k}(y^{k-1}) = x^k$ и формулу (2), получим

$$\sum_{k=1}^{N} \gamma_k (x^k - x)^T \nabla_x f_k(x^k; \xi^k) \leq W_{\beta_1}(y^0) + x^T y^N - W_{\beta_{N+1}}(y^N) + \sum_{k=1}^{N} \frac{\gamma_k^2}{2\alpha\beta_k} \|\nabla_x f_k(x^k; \xi^k)\|_\infty^2,$$

где $x^T$ – означает транспонирование вектора $x$, который мы выбираем так, чтобы он доставлял решение задачи (1). Поскольку [4]

$$W_{\beta_1}(y^0) = W_{\beta_1}(0) = 0 \text{ и } \beta_{N+1} V(x) \geq x^T y^N - W_{\beta_{N+1}}(y^N),$$

то

$$\sum_{k=1}^{N} \gamma_k (x^k - x)^T \nabla_x f_k(x^k; \xi^k) \leq \beta_{N+1} V(x) + \sum_{k=1}^{N} \frac{\gamma_k^2}{2\alpha\beta_k} \|\nabla_x f_k(x^k; \xi^k)\|_\infty^2.$$

Тогда, из выпуклости функции $E_{\xi^k}[f_k(x; \xi^k)]$ по $x$ (в виду условия 1) следует, что

$$\sum_{k=1}^{N} \gamma_k \left\{ E_{\xi^k}[f_k(x^k; \xi^k)] - E_{\xi^k}[f_k(x; \xi^k)] \right\} \leq$$
$$\leq \sum_{k=1}^{N} \gamma_k (x^k - x)^T \nabla_x E_{\xi^k}[f_k(x^k; \xi^k)] \leq \qquad (3)$$
$$\leq \beta_{N+1} V(x) - \sum_{k=1}^{N} \gamma_k (x^k - x)^T \left( \nabla_x f_k(x^k; \xi^k) - \nabla_x E_{\xi^k}[f_k(x^k; \xi^k)] \right) + \sum_{k=1}^{N} \frac{\gamma_k^2}{2\alpha\beta_k} \|\nabla_x f_k(x^k; \xi^k)\|_\infty^2.$$

Возьмем полное (т.е., в отличие от замечания 2, с учетом зависимости $x^k$ от $\{\xi^1,...,\xi^{k-1}\}$) математическое ожидание (в два шага $E[\cdot] = E\left[E_{\xi^k}[\cdot | \Xi^{k-1}]\right]$ – для каждого слагаемого свое $k$) от обеих частей неравенства, учитывая условие 3.а и то, что

$$E\left((x^k - x)^T \left(\nabla_x f_k(x^k; \xi^k) - \nabla_x E_{\xi^k}[f_k(x^k; \xi^k)]\right)\right) =$$
$$= E\left[(x^k - x)^T E_{\xi^k}\left(\nabla_x f_k(x^k; \xi^k) - \nabla_x E_{\xi^k}[f_k(x^k; \xi^k)] \big| \Xi^{k-1}\right)\right] = 0,$$

поскольку $x^k$ – $\Xi^{k-1}$-измеримый вектор и внутреннее условное математическое ожидание в силу условия 2 равно 0, получим

$$\sum_{k=1}^{N} \gamma_k \left\{ E\left[ f_k\left(x^k; \xi^k\right)\right] - E_{\xi^k}\left[ f_k\left(x; \xi^k\right)\right]\right\} \leq \beta_{N+1} V(x) + \sum_{k=1}^{N} \frac{\gamma_k^2}{2\alpha\beta_k} E\left[\left\|\nabla_x f_k\left(x^k; \xi^k\right)\right\|_\infty^2\right] \leq$$

$$\leq \beta_{N+1} R^2 + M^2 \sum_{k=1}^{N} \frac{\gamma_k^2}{2\alpha\beta_k},$$

где $R^2 = \max_{x \in S_n(1)} V(x) \equiv \ln n$. Подставляя $\gamma_k \equiv 1$, и минимизируя правую часть неравенства по неубывающим последовательностям с положительными элементами $\{\beta_k\}_{k=1}^{N+1}$, не допуская при этом зависимость $\{\beta_k\}_{k=1}^{N+1}$ от потенциально неизвестного $N$, получим

$$\beta_k = \sqrt{\frac{M^2}{\alpha R^2}} \sqrt{k},$$

$$\frac{1}{N} \sum_{k=1}^{N} E\left[ f_k\left(x^k; \xi^k\right)\right] - \min_{x \in S_n(1)} \frac{1}{N} \sum_{k=1}^{N} E_{\xi^k}\left[ f_k\left(x; \xi^k\right)\right] \leq 2\sqrt{\frac{M^2 R^2}{\alpha N}}.$$

Для того чтобы доказать первую часть теоремы, осталось подставить $\alpha = 1$ и $R^2 = \ln n$.

**Замечание 3.** Если разрешать $\{\beta_k\}_{k=1}^{N+1}$ зависеть от $N$ (см., например, алгоритм МЗС2-неадаптивный в следующем пункте), то в последней формуле "2"-у можно занести под знак корня [4]. Все это переносится и на последующие рассуждения с вероятностями больших отклонений.

**Замечание 4.** Строго говоря, в получено оценке в знаменателе вместо $N$ нужно писать $N^2/(N+1)$. Считая, что $N \gg 1$, мы пренебрегли этим для компактности записи.

Для доказательства второй части теоремы вернемся к формуле (3). Из условия 3 имеем

$$P\left( \sum_{k=1}^{N} \frac{\gamma_k^2}{2\alpha\beta_k} \left\|\nabla_x f_k\left(x^k; \xi^k\right)\right\|_\infty^2 > M^2 \sum_{k=1}^{N} \frac{\gamma_k^2}{2\alpha\beta_k} \right) = 0. \quad (4)$$

Из неравенства Азума–Хефдинга [21], подобно [7], [22], получаем для ограниченной мартингал-разности

$$\left|\left(x - x^k\right)^T \left(\nabla_x f_k\left(x^k; \xi^k\right) - \nabla_x E_{\xi^k}\left[ f_k\left(x^k; \xi^k\right)\right]\right)\right| \leq 4M$$

следующее неравенство:

$$P\left( \sum_{k=1}^{N} \gamma_k \left(x - x^k\right)^T \left(\nabla_x f_k\left(x^k; \xi^k\right) - \nabla_x E_{\xi^k}\left[ f_k\left(x^k; \xi^k\right)\right]\right) \geq 4M\Lambda \sqrt{\sum_{k=1}^{N} \gamma_k^2} \right) \leq \exp\left(-\Lambda^2/2\right).$$

Подставляя $\gamma_k \equiv 1$, $\Lambda = \sqrt{2\Omega}$, получим вторую часть теоремы.

**Замечание 5.** В статье [7] для задачи стохастической выпуклой оптимизации ($\{\xi^k\}$ – независимые случайные величины) приводится оценка вероятностей больших уклонений с точностью до констант аналогичная оценке, приведенной в теореме. Аналогичная оценка приводится в [7] и для случая, когда вместо условия 3 предполагается условие 3.б. При этом в [7] использовалось условие независимости $\{\xi^k\}$ (при установлении неравенства типа (4) в общем случае и в неравенстве Азума–Хефдинга). Для онлайн оптимизации, как правило, независимость $\{\xi^k\}$ место не имеет. Тем не менее, условия 2, 3.б обеспечивают выполнения этих неравенств в том же виде, как если бы независимость $\{\xi^k\}$ имела место. Приведем теперь оценки в случае тяжелых хвостов. Если $\|\nabla_x f(x,\xi)\|_\infty^2$ имеет степенной хвост ($\alpha > 2$)

$$P\left(\frac{\|\nabla_x f(x,\xi)\|_\infty^2}{M^2} \geq t\right) = \mathrm{O}\left(\frac{1}{t^\alpha}\right),$$

то существует такая константа $C_\alpha > 0$, что с вероятностью $\geq 1-\sigma$

$$\frac{1}{N}\sum_{k=1}^N E_{\xi^k}\left[f_k\left(x^k;\xi^k\right)\right] - \min_{x \in S_n(1)} \frac{1}{N}\sum_{k=1}^N E_{\xi^k}\left[f_k\left(x;\xi^k\right)\right] \leq C_\alpha M \frac{\sqrt{\ln n \ln(\sigma^{-1})} + \frac{(N/\sigma)^{1/\alpha}}{N}}{\sqrt{N}}.$$

Если мы не делаем никаких предположений относительно распределений случайных величин $\{\xi^k\}$, кроме 1, 2 и существования первых двух равномерно ограниченных моментов у $\nabla_x f_k(x^k;\xi^k)$, то из неравенства Маркова и первого неравенства в теореме 1 (на математические ожидания) имеем: существует такая константа $C > 0$, что с вероятностью $\geq 1-\sigma$

$$\frac{1}{N}\sum_{k=1}^N E_{\xi^k}\left[f_k\left(x^k;\xi^k\right)\right] - \min_{x \in S_n(1)} \frac{1}{N}\sum_{k=1}^N E_{\xi^k}\left[f_k\left(x;\xi^k\right)\right] \leq \frac{CM}{\sigma}\sqrt{\frac{\ln n}{N}}.$$

Труднее обстоит дело, если мы хотим оценить регрет или

$$\frac{1}{N}\sum_{k=1}^N f_k\left(x^k;\xi^k\right) - \min_{x \in S_n(1)} \frac{1}{N}\sum_{k=1}^N f_k\left(x;\xi^k\right).$$

Тем не менее, при дополнительных оговорках и такие выражения можно вероятностно оценивать [14], [22].

Аналогичное замечание имеет место и для теоремы 2 ниже.

## 3. Онлайн МЗС со стохастической проекцией

Снова рассмотрим постановку задачи стохастической онлайн оптимизации (1) из п. 2. Но на этот раз будем допускать, что метод генерирования последовательности $\{x^k\}$ может допускать (внешнюю, дополнительную) рандомизацию. Это допущение позволит частично перенести результаты п. 2 на не выпуклые функции $E_{\xi^k}\left[f_k\left(x;\xi^k\right)\right]$ (см. пример 4 п. 4), на ситуации, когда по условию задачи $\{x^k\}$ должны выбираться среди вершин единичного симплекса (примеры 1 и 5 п. 4). Также как и раньше онлайновость постановки задачи допускает, что на каждом шаге $k$ функция $f_k$ может подбираться из рассматриваемого класса функций враждебно по отношению к используемому нами методу генерации последовательности $\{x^k\}$. В частности, $f_k$ может зависеть от $\{x^1, \xi^1, f_1(\cdot);...;x^{k-1}, \xi^{k-1}, f_{k-1}(\cdot)\}$, и даже от распределения вероятностей $p^k$ (многорукие бандиты), согласно которому осуществляется выбор $x^k$. Чтобы можно было работать с таким классом задач, нам придется наложить дополнительное **условие**:

4. На каждом шаге генерирование случайной величины $x^k$ согласно распределению вероятностей $p^k$ осуществляется независимо ни от чего. Выбор $f_k$ осуществляется без знания реализации $x^k$.

Положим $p_i^1 = x_i^1 = 1/n$, $i = 1,...,n$. Пусть $t = 1,...,N-1$.

### Алгоритм МЗС2-адаптивный / Метод двойственных усреднений

*Согласно распределению вероятностей*

$$p_i^{t+1} = \frac{\exp\left(-\dfrac{1}{\beta_{t+1}}\sum_{k=1}^{t}\dfrac{\partial f_k\left(x^k;\xi^k\right)}{\partial x_i}\right)}{\sum_{l=1}^{n}\exp\left(-\dfrac{1}{\beta_{t+1}}\sum_{k=1}^{t}\dfrac{\partial f_k\left(x^k;\xi^k\right)}{\partial x_l}\right)}, \quad i = 1,...,n, \, \beta_t = \frac{M\sqrt{t}}{\sqrt{\ln n}},$$

*получаем случайную величину* $i(t+1)$, $x_{i(t+1)}^{t+1} = 1$, $x_j^{t+1} = 0$, $j \neq i(t+1)$.

### Алгоритм МЗС2-неадаптивный (заранее известно $N$)

*Согласно распределению вероятностей*

$$p_i^{t+1} = \frac{\exp\left(-\dfrac{1}{\beta_{t+1}}\sum_{k=1}^{t}\gamma_k\dfrac{\partial f_k\left(x^k;\xi^k\right)}{\partial x_i}\right)}{\sum_{l=1}^{n}\exp\left(-\dfrac{1}{\beta_{t+1}}\sum_{k=1}^{t}\gamma_k\dfrac{\partial f_k\left(x^k;\xi^k\right)}{\partial x_l}\right)}, \quad i = 1,...,n,$$

$$\gamma_k \equiv M^{-1}\sqrt{2\ln n/N}, \ \beta_t \equiv 1,$$

*получаем случайную величину* $i(t+1)$, $x^{t+1}_{i(t+1)} = 1$, $x^{t+1}_j = 0$, $j \neq i(t+1)$.

**Мотивация (ограничимся детерминированным случаем с $\gamma_k \equiv 1$).**
Аппроксимируя

$$\min_{x \in S_n(1)} \frac{1}{N}\sum_{k=1}^{t} f_k(x) \approx \min_{x \in S_n(1)} \frac{1}{N}\sum_{k=1}^{t}\left\{ f_k(x^k) + \langle \nabla f_k(x^k), x - x^k \rangle \right\},$$

получим $x^{t+1}_j = 1; x^{t+1}_i = 0, i \neq j$ (для простоты считаем, что имеет место не вырожденный случай)

$$j = \arg\max_{i=1,\ldots,n}\left\{ \sum_{k=1}^{t}\left[-\nabla f_k(x^k)\right]_i \right\}.$$

Поскольку мы работаем с аппроксимацией (нижними аффинными минорантами) исходной задачи, то предлагается немного видоизменить это правило

$$P\left(x^{t+1}_j = 1; x^{t+1}_i = 0, i \neq j\right) \stackrel{def}{=} P_\varsigma\left( j = \arg\max_{i=1,\ldots,n}\left\{\left(\sum_{k=1}^{t}\left[-\nabla f_k(x^k)\right]_i\right) + \varsigma_{t,i}\right\}\right),$$

где $\varsigma_{t,i}$ – независимые одинаково распределенные случайные величины по закону Гумбеля [23], [24] с параметром $\beta_{t+1}$, характеризующим среднеквадратичное отклонение $\varsigma_{t,i}$:

$$P(\varsigma_{t,i} < \tau) = \exp\left\{-e^{-\tau/\beta_{t+1}}\right\}.$$

Можно показать [3], [9], [23], что

$$E_\varsigma\left[x^{t+1}\right] = \nabla W_{\beta_{t+1}}\left(-\sum_{k=1}^{t}\nabla f_k(x^k)\right).$$

**Замечание 6.** Естественно задаться вопросом: а какое распределение "наиболее подходит" для $\varsigma_{t,i}$, чтобы в случае "враждебной Природы" (то есть в минимаксном смысле) иметь наилучшие оценки [9]? Ответом будет [9]: показательное распределение (точнее $\beta_{t+1}\text{Exp}(1)$, где $\text{Exp}(1)$ – случайная величина, имеющая показательное распределение с параметром 1), которое ведет себя в интересном для анализа диапазоне подобно распределению Гумбеля, но в случае Гумбеля мы явно можем посчитать интересующие нас вероятности. Как правило, такого рода задачи явно не решаются, и распределение Гумбеля является приятным исключением, для которого есть явные формулы.

Приведенная выше мотивация имеет одно интересное приложение в содержательной интерпретации равновесного распределения транспортных потоков. Не много об этом написано [25].

К сожалению, не делая относительно функций $f_k(x;\xi^k)$ дополнительно никаких предположений, не удается доказать для МЗС2 аналог теоремы 1. Чтобы можно было сформулировать такой аналог, мы вынуждены будем предполагать, что $f_k(x;\xi^k)$ – линейные функции по $x$ (можно обобщить и на сублинейные). С одной стороны это существенно сужает класс задач, к которым применим МЗС2. С другой стороны, как будет продемонстрировано в следующем пункте, даже такой узкий класс функций за счет "онлайновости" позволяет применять МЗС2 к довольно широкому кругу задач. Для того чтобы лучше чувствовалась преемственность методов и доказательств их сходимости, далее мы по-прежнему будем использовать общие обозначения $f_k(x;\xi^k)$, не подчеркивая в формулах линейность.

**Теорема 2.** *Пусть справедливы условия 1, 2, 3.а, 4 и $f_k(x;\xi^k)$ – линейные функции по $x$, тогда*

$$\frac{1}{N}\sum_{k=1}^{N} E\left[f_k\left(x^k;\xi^k\right)\right] - \min_{x\in S_n(1)}\frac{1}{N}\sum_{k=1}^{N} E_{\xi^k}\left[f_k\left(x;\xi^k\right)\right] \le 2M\sqrt{\frac{\ln n}{N}}.$$

*Для неадаптивного метода "2"-у перед $M$ можно занести под знак корня.*

*Кроме того, если справедливы условия 1, 2, 3, 4, то при $\Omega \ge 0$*

$$P_{x^1,...,x^N}\left\{\frac{1}{N}\sum_{k=1}^{N} E_{\xi^k}\left[f_k\left(x^k;\xi^k\right)\right] - \min_{x\in S_n(1)}\frac{1}{N}\sum_{k=1}^{N} E_{\xi^k}\left[f_k\left(x;\xi^k\right)\right] \ge \frac{2M}{\sqrt{N}}\left(\sqrt{\ln n}+\sqrt{18\Omega}\right)\right\} \le \exp(-\Omega).$$

*Если $f_k(x;\xi^k) \equiv f_k(x)$, то это неравенство можно уточнить*

$$\frac{2M}{\sqrt{N}}\left(\sqrt{\ln n}+\sqrt{18\Omega}\right) \to \frac{2M}{\sqrt{N}}\left(\sqrt{\ln n}+\sqrt{2\Omega}\right),$$

*при этом же условии для неадаптивного метода можно еще больше уточнить*

$$\frac{2M}{\sqrt{N}}\left(\sqrt{\ln n}+\sqrt{2\Omega}\right) \to \frac{\sqrt{2}M}{\sqrt{N}}\left(\sqrt{\ln n}+2\sqrt{\Omega}\right).$$

**Схема доказательства теоремы 2.** Доказательство фактически дословно повторяет доказательство теоремы 1. Небольшая разница лишь в том, что основная формула (3) перепишется следующим образом:

$$\sum_{k=1}^{N} \gamma_k \left\{ E_{\xi^k}\left[ f_k\left(x^k;\xi^k\right)\right] - E_{\xi^k}\left[ f_k\left(x;\xi^k\right)\right]\right\} \le \sum_{k=1}^{N} \gamma_k \left(x^k - x\right)^T \nabla_x E_{\xi^k}\left[ f_k\left(x^k;\xi^k\right)\right] \le$$

$$\le \beta_{N+1} V(x) + \sum_{k=1}^{N} \gamma_k \left(x^k - E_{x^k}\left[x^k\right]\right)^T \nabla_x f_k\left(x^k;\xi^k\right) -$$

$$- \sum_{k=1}^{N} \gamma_k \left(x^k - x\right)^T \left(\nabla_x f_k\left(x^k;\xi^k\right) - \nabla_x E_{\xi^k}\left[ f_k\left(x^k;\xi^k\right)\right]\right) + \sum_{k=1}^{N} \frac{\gamma_k^2}{2\alpha\beta_k} \left\| \nabla_x f_k\left(x^k;\xi^k\right)\right\|_\infty^2.$$

Здесь мы просто не много по-другому (по сравнению с доказательством теоремы 1) переписали неравенство

$$\sum_{k=1}^{N} \gamma_k \left(x^k - x\right)^T \nabla_x f_k\left(x^k;\xi^k\right) \le \beta_{N+1} V(x) + \sum_{k=1}^{N} \frac{\gamma_k^2}{2\alpha\beta_k} \left\| \nabla_x f_k\left(x^k;\xi^k\right)\right\|_\infty^2,$$

используя выпуклость функции $E_{\xi^k}\left[ f_k\left(x;\xi^k\right)\right]$ по $x$ в виду условия 1.

Введем случайные величины

$$Z_k = \left(x^k - E_{x^k}\left[x^k\right]\right)^T \nabla_x f_k\left(x^k;\xi^k\right), \quad \tilde{Z}_k = \left(x^k - x\right)^T \left(\nabla_x f_k\left(x^k;\xi^k\right) - \nabla_x E_{\xi^k}\left[ f_k\left(x^k;\xi^k\right)\right]\right)$$

Поскольку $f_k\left(x^k;\xi^k\right) = l_k\left(\xi^k\right)^T x^k$ – линейные функции, то по условиям 4, 2

$$E_{x^k}\left[ Z_k | \Xi^{k-1}\right] \equiv 0, \quad E_{\xi^k}\left[ \tilde{Z}_k | \Xi^{k-1}\right] \equiv 0.$$

Именно в этом месте и только в нем используется линейность $f_k\left(x;\xi^k\right)$. К сожалению, предложенный здесь способ рассуждения не позволяет хоть сколько-нибудь ослабить это условие.

Рассуждая далее также как в доказательстве теоремы 1, получим теорему 2.

**Замечание 7.** Как уже отмечалось во введении, теорема 2 во многом мотивирована работой [6]. Собственно, форма, в который мы представили алгоритм МЗС2-неадаптивный, выбрана именно такой (альтернативным вариантом было положить $\gamma_k \equiv 1$, $\beta_k \equiv M\sqrt{N/(2\ln n)}$), чтобы была максимальная близость к алгоритму работы [6].

**Замечание 8.** Идея рандомизации (искусственного введения случайности), положенная в основу описанных алгоритмов, чрезвычайно продуктивна: против нас играет, возможно, враждебная "Природа", которая, зная историю игры, и наши текущие намерения старается нам "предложить вариант похуже". С этим можно "бороться" за счет случайного независимого осуществления своих намерений на каждом шаге, с реализацией неизвестной "Природе". За счет этой случайности мы переходим от анализа по худшему случаю (роль которого в онлайн оптимизации играет враждебная "Природа") к анализу "в среднем". Такая рандомизация, как будет отмечено в следующем пункте, дает возможность получать оценки, которые в детерминированном случае получить невозможно. Причем, если в онлайн постановке такая рандомизация прописывается в "правилах игры", то применительно к задачам обычной оптимизации все это возникает совершенно естественным образом, как желание с большой вероятностью обезопасить себя от "самых худших случаев" детерминированной версии метода. Отметим, что речь идет о, так называемых, массовых задачах, т.е., исследуя тот или иной метод, мы точно не знаем какой конкретно объект поступит на вход, поэтому, чтобы гарантировано что-то иметь, мы исходим из худшего (наименее благоприятного для данного метода) случая входных данных. Описанный МЗС2 естественно также понимать как покомпонентный метод (стохастического) субградиентного спуска со случайным выбором компоненты. Происходит рандомизация при проектировании (в смысле расстояния Брэгмана) на единичный симплекс. А именно, если проектироваться в указанном выше смысле на единичный симплекс, то получится вектор, который можно проинтерпретировать как распределение вероятностей некоторой дискретной случайной величины, принимающей значения $1,...,n$. Если выбрать вершину симплекса согласно этой дискретной случайной величине, и заменить проекцию этой вершиной, то получим случайную проекцию, математическое ожидание которой равно честной проекции. Как будет отмечено в следующем пункте, такой метод не только оптимален с точки зрения числа итераций, но и в некотором смысле с точки зрения затрат на выполнение одной итерации (см. также [26]).

**Замечание 9 (Ф.А. Федоренко).** Если ограничение $x \in S_n(1)$ в задаче (1) заменить на ограничение

$$x = \left(z^1,...,z^m\right) \in \prod_{j=1}^{m} S_{n_j}\left(d_j\right),$$

типичное для популяционных игр (с несколькими популяциями), в частности, для транспортных приложений [25], то взяв

$$V(x) = \sum_{j=1}^{m} V_j\left(z^j\right), \; V_j\left(z^j\right) = d_j \ln n_j + \sum_{i=1}^{n_j} z_i^j \ln\left(\frac{z_i^j}{d_j}\right), \; W_\beta(y) = \sup_{z^j \in S_{n_j}(d_j)} \left\{ \sum_{j=1}^{m} \left\langle v^j, z^j \right\rangle - \beta V(x) \right\},$$

где $y = \left(v^1,...,v^m\right)$, можно получить аналогичные теореме 1 оценки с заменой (аналогично по теореме 2)

$$2M\sqrt{\frac{\ln n}{N}} \to 2M\sqrt{\frac{\left(\max_{j=1,...,m} d_j\right) \cdot \left(\sum_{j=1}^{m} d_j \ln n_j\right)}{N}},$$

$$\frac{2M}{\sqrt{N}}\left(\sqrt{\ln n}+\sqrt{8\Omega}\right) \to \frac{2M}{\sqrt{N}}\left(\sqrt{\left(\max_{j=1,...,m} d_j\right)\cdot\left(\sum_{j=1}^{m} d_j \ln n_j\right)} + \left(\sum_{j=1}^{m} d_j\right)\cdot\sqrt{8\Omega}\right).$$

## 4. Приложения МЗС

В заключительном пункте мы постараемся продемонстрировать некоторые возможности и ограничения описанных в пп. 2 и 3 методов. Мы не будем стремиться здесь к максимальной общности или рассмотрению всех основных приложений МЗС.

**Пример 1 (многорукие бандиты [5], [9], [14]).** Имеется $n$ различных ручек. Игра повторяется $N \gg 1$ раз (это число может быть заранее неизвестно). На каждом шаге $k$ мы должны выбрать ручку $i(k)$, которую "дергаем". Дергание ручки приносит нам некоторые, вообще говоря, случайные потери $r^k_{i(k)}$ (считаем, для определенности, что всегда $r^k_{i(k)} \in [0,1]$), зависящие от номера шага, номера ручки и от того, какой стратегии мы придерживались до шага $k$ включительно. Наша стратегия на шаге $k$ описывается вектором распределения вероятностей $x^k \in S_n(1)$, согласно которому мы независимо ни от чего выбираем ручку, которую будем дергать. Все, чем мы располагаем на шаге $k$, это вектором

$$\left(\left(x^1, i(1), r^1_{i(1)}\right); ...; \left(x^{k-1}, i(k-1), r^{k-1}_{i(k-1)}\right)\right).$$

Мы считаем, что потери на $k$-м шаге $r^k$ зависят от $x^k$ (но не от результата разыгрывания из распределения $x^k$), зависят от $(x^1,...,x^{k-1})$ и результатов соответствующих разыгрываний, а также зависят от $(r^1,...,r^{k-1})$. Целью является таким образом организовать процедуру дергания ручек, чтобы ожидаемые суммарные потери были бы минимальны. Введем функцию ($r^k$ и результат разыгрывания, согласно распределению вероятностей, заданному вектором $x$, – независимы; обе эти "случайности" мы обозначаем $\xi^k$)

$$f_k(x; \xi^k) = r^k_i \text{ с вероятностью } x_i, \ i=1,...,n,$$

и её обобщенный (в смысле удовлетворения условию 2) стохастический градиент

$$\nabla_x f_k(x; \xi^k) = (0,...,\underbrace{r^k_i}_{i}/x_i,...,0)^T \text{ с вероятностью } x_i, \ i=1,...,n.$$

Тогда выполнены условия 1 и 2. Однако имеется проблема: константа $M$ в условии 3 получается слишком большой (например, в 3.а $M = \sup_{x \in S_n(1)} \sqrt{\sum_{i=1}^{n} x_i^{-1}} = \infty$), то есть теорема 1 ничего дать не может в том виде, в котором она была нами приведена. Возникает желание "что-то подкрутить" в доказательстве теоремы, чтобы можно было ей воспользоваться. Уже по ходу самого доказательства мы отмечали неравенство (*), которое может оказаться довольно грубым в определенных ситуациях. Многорукие бандиты дают пример как раз такой ситуации. Более аккуратный анализ [5], [10], [14], использующий специфику данной задачи, позволяет оценить (см. доказательство теоремы 1)

$$\gamma_k \nabla_x f_k\left(x^k; \xi^k\right)^T \int_0^1 \left(\nabla W_{\beta_k}\left(\tau y^k + (1-\tau) y^{k-1}\right) - \nabla W_{\beta_k}\left(y^{k-1}\right)\right) d\tau$$

точнее, что приводит в основной формуле (3) к замене слагаемых вида

$$\frac{\gamma_k^2}{2\alpha\beta_k} \left\|\nabla_x f_k\left(x^k; \xi^k\right)\right\|_\infty^2$$

на (в этом месте, для наглядности, мы намеренно несколько упрощаем и огрубляем)

$$\gamma_k^2 \frac{x_j^k\left(1-x_j^k\right)}{\alpha\beta_k}\left(\frac{r_j^k}{x_j^k}\right)^2,$$

где $j$ – номер ручки, выбранной алгоритмом на $k$-м шаге. В результате мы получаем, что теорема 1 остается верной с эффективной константой $M = \sqrt{2n}$. Таким образом, действуя согласно МЗС1, наши потери (псевдо регрет [14]) будут

$$O\left(\sqrt{\frac{n\ln n}{N}}\right) - \text{в среднем;} \quad O\left(\sqrt{\frac{n\ln(n/\sigma)}{N}}\right) - \text{с вероятностью} \geq 1-\sigma,$$

что с точностью до логарифмического фактора соответствует нижним оценкам [5], [9], [14].

**Замечание 10.** Стоит обратить внимание, что если использовать более специальную прокс-структуру [5], [14], то для псевдо регрета можно получить оценки без логарифмического фактора $\ln n$ под корнем, что уже соответствует нижним оценкам. В частности, это обстоятельство означает, что выбирать "расстояние" Брэгмана для симплекса не всегда оптимально (но близко к оптимуму). Тем не менее, в последующих нескольких примерах мы убедимся, что для ряда других постановок, оценки, полученные с помощью прокс-структуры, порожденной расстоянием Брэгмана, – оптимальные. Кроме того, есть еще плата за избавление от фактора $\ln n$ под корнем – удорожание процедуры вычисления проекции на симплекс в смысле этой прокс-структуры. Другими словами, это ускорение оправдано только для онлайн постановок, в которых, как правило, стремятся минимизировать (псевдо) регрет, не сильно учитывая общую вычислительную трудоемкость.

Труднее обстоит дело с оценкой регрета (в среднем и вероятностей больших уклонений). Пример из лекции 6 [11] показывает ($n = 2$), что МЗС1 может давать регрет $\sim cN^{-1/4}$, что значительно хуже оценки псевдо регрета $\sim cN^{-1/2}$. По сути, речь идет о том, что написано в конце замечания 5. Здесь уже требуется некая игра "bias–variance trade off": отказаться от несмещенности оценки градиента для уменьшения дисперсии этой оценки. Этот популярный трюк в математической статистике и машинном обучении позволяет с некоторыми оговорками распространить приведенные выше оценки псевдо регрета $O\left(\sqrt{n\ln(n/\sigma)/N}\right)$ и на случай оценок регрета. Кое-что на эту тему применительно к многоруким бандитам можно найти в обзоре [14].

Интересно заметить, что результаты, описанные в примере 1 можно получить (с аналогичными оговорками) с помощью МЗС2 и теоремы 2. Для этого нужно взять

$$f_k(x;\xi^k) = \langle r^k, x \rangle$$

и её обобщенный (в смысле удовлетворения условию 2) стохастический градиент

$$\nabla_x f_k(x;\xi^k) = (0,...,\underbrace{r_i^k}_{i}/p_i,...,0)^T, \text{ если } x = (0,...,\underbrace{1}_{i},...,0)^T,$$

где $x = (0,...,\underbrace{1}_{i},...,0)^T$ с вероятностью $p_i$, $i=1,...,n$,

здесь $\xi^k$ отражает только случайность, сидящую в $r^k$. Это определение стохастического градиента (в отличие от рассмотренного выше) явно учитывает распределение вероятностей $p$, из которого генерируется номер единственной ненулевой (единичной) компоненты вектора $x$. Также как и раньше для выполнения условия 2 необходимо предполагать независимость $\xi^k$ и процедуры разыгрывания единичный компоненты вектора $x$ согласно распределению $p$.

**Пример 2 (взвешивание экспертных решений, линейные потери [9], [27]).** Рассмотрим задачу взвешивание экспертных решений, следуя [9], [27]. Имеется $n$ различных Экспертов. Каждый Эксперт играет на рынке. Игра повторяется $N \gg 1$ раз (это число может быть заранее неизвестно). Пусть $l_i^k$ – проигрыш Эксперта $i$ на шаге $k$ ($\left|l_i^k\right| \le M$). На каждом шаге $k$ мы распределяем один доллар между Экспертами, согласно вектору $x^k \in S_n(1)$. Потери, которые мы при этом несем, рассчитываются по потерям экспертов $\langle l^k, x^k \rangle$. Целью является таким образом организовать процедуру распределения доллара на каждом шаге, чтобы наши суммарные потери были бы минимальны. Допускается, что потери экспертов $l^k$ могут зависеть еще и от текущего хода $x^k$. Легко проверить, что для данной постановки применима теорема 1 в детерминированном варианте с функциями

$$f_k(x;\xi^k) \equiv f_k(x) = \langle l^k, x \rangle.$$

При этом оценка, даваемая теоремой 1,

$$\mathrm{O}\left(M\sqrt{\frac{\ln n}{N}}\right)$$

– оптимальна для данного класса задач [9], [27].

**Пример 3 (взвешивание экспертных решений, выпуклые потери [9], [27]).** В условиях предыдущего примера предположим, что на $k$-м шаге $i$-й эксперт использует стратегию $\zeta_i^k \in \Delta$ (множество $\Delta$ – выпуклое), дающую потери $\lambda\left(\omega^k, \zeta_i^k\right)$, где $\omega^k$ – "ход", возможно, враждебной "Природы", знающей, в том числе, и нашу текущую стратегию. Функция $\lambda(\cdot)$ – выпуклая по второму аргументу и $|\lambda(\cdot)| \le M$. На каждом шаге мы должны выбирать свою стратегию

$$x \stackrel{def}{=} \sum_{i=1}^{n} x_i \cdot \zeta_i^k \in \Delta,$$

дающую потери $\lambda\left(\omega^k, x\right)$ так, чтобы наши суммарные потери были минимальны. Для данной постановки также применима теорема в детерминированном варианте с

$$f_k\left(x; \xi^k\right) \equiv f_k(x) = \sum_{i=1}^{n} x_i \lambda\left(\omega^k, \zeta_i^k\right) \ge \lambda\left(\omega^k, x\right).$$

Чтобы применить теорему осталось заметить, что функция $\lambda\left(\omega^k, \zeta\right)$ – выпуклая по $\zeta$ для любого $\omega^k$, поэтому

$$\sum_{k=1}^{N} \lambda\left(\omega^k, x^k\right) - \min_{i=1,\ldots,n} \sum_{k=1}^{N} \lambda\left(\omega^k, \zeta_i^k\right) \le \sum_{k=1}^{N} f_k\left(x^k\right) - \min_{x \in S_n(1)} \sum_{k=1}^{N} f_k(x).$$

При этом оценка, даваемая теоремой 1,

$$\mathrm{O}\left(M\sqrt{\frac{\ln n}{N}}\right)$$

– также оптимальна для данного класса задач [9], [27].

Полезно, на наш взгляд, будет здесь привести другой способ (более типичный для данного класса приложений) получения аналогичного результата, не связанный на прямую с МЗС схемой вывода, но фактически, приводящий к точно такому же алгоритму. Этот способ также весьма популярен в машинном обучении, теории игр, теории алгоритмов [9], [27], [28].[2]

---

[2] Максимум из независимых случайных величин, который сложно исследовать, заменяется (с хорошей точностью, контролируемой малостью параметра $\beta$) логарифмом от суммы экспонент от этих независимых случайных величин. А сумму независимых случайных величин (их экспонент) исследовать уже на много проще. В оптимизации эту процедуру называют сглаживанием [29].

Введем обозначение $L_i^N = \sum_{k=1}^N \lambda(\omega^k, \zeta_i^k)$, $\tilde{L}^N = \sum_{k=1}^N \lambda(\omega^k, x^k)$, по определению считаем $L_i^0 \equiv 0$. Рассмотрим

$$W_\beta\left(\{-L_i^N\}_{i=1}^n\right) = \beta \ln\left(\frac{1}{n}\sum_{i=1}^n \exp(-L_i^N/\beta)\right) \geq -\min_{i=1,\ldots,n} L_i^N - \beta \ln n.$$

С другой стороны, вводя дискретную случайную величину $z^k$, имеющую (независящее ни от чего) распределение $x^k$ (рассчитанное также как и раньше исходя из МЗС1, примененного к набору функций $\{f_k(x)\}_{k=1}^N$, определенных выше), можно заметить, что

$$W_\beta\left(\{-L_i^N\}_{i=1}^n\right) = \sum_{k=1}^N \left(W_\beta\left(\{-L_i^k\}_{i=1}^n\right) - W_\beta\left(\{-L_i^{k-1}\}_{i=1}^n\right)\right) = \beta \sum_{k=1}^N \ln\left(E_z\left(e^{-\lambda(\omega^k, z^k)/\beta}\right)\right).$$

Используя далее неравенство Хефдинга (для с.в. $X \in [-M, M]$), см. [21]

$$\ln\left(E_X\left(e^{sX}\right)\right) \leq sE_X(X) + s^2 \frac{M^2}{2},$$

Получим

$$W_\beta\left(\{-L_i^N\}_{i=1}^n\right) \leq -\tilde{L}^N + (2\beta)^{-1} M^2 N.$$

Таким образом,

$$\tilde{L}^N \leq \min_{i=1,\ldots,n} L_i^N + \beta \ln n + (2\beta)^{-1} M^2 N.$$

Минимизация правой части по $\beta > 0$ приводит нас к уже известному ответу. Аналогичные, но чуть более тонкие рассуждения, позволяют избавиться от зависимости $\beta$ от $N$, то есть сделать алгоритм адаптивным.

**Пример 4 (взвешивание экспертных решений, невыпуклые потери [9], [27]).** Предположим, что в условиях примера 3 мы не можем гарантировать выпуклость $\lambda(\cdot)$ – по второму аргументу. Тогда мы выбираем стратегию – распределение вероятностей на множестве стратегий Экспертов, и разыгрываем случайную величину согласно этому распределению вероятностей. Другими словами мы просто пользуемся МЗС2 с $f_k(x; \xi^k) \equiv f_k(x) = \sum_{i=1}^n x_i \lambda(\omega^k, \zeta_i^k)$, применимость которого обосновывается теоремой 2, с оценками

$$\mathrm{O}\left(M\sqrt{\frac{\ln n}{N}}\right) \text{ – в среднем; } \mathrm{O}\left(M\sqrt{\frac{\ln(n/\sigma)}{N}}\right) \text{ – с вероятностью } \geq 1-\sigma,$$

– оптимальными для данного класса задач [9], [27]. Ключевая разница в примерах 1 и 4, "сто́ящая" $\sim \sqrt{n}$ в оценке $M$ для многоруких бандитов (пример 1), заключается в том, что в многоруких бандитах мы имеем только свою историю дергания ручек (нам не известно, какие бы потери нам принесли другие ручки, кабы мы их выбрали), а в постановке взвешивания экспертных решений это все известно, и называется потерями экспертов.

Как будет видно из следующего примера описанный только что подход вполне успешно работает (дает не улучшаемые результаты) и в случае выпуклой по второму аргументу функции $\lambda(\cdot)$.

**Пример 5 (антагонистические матричные игры [6], [9], [26], [27]).** Пусть есть два игрока А и Б. Задана матрица игры $A = \|a_{ij}\|$, где $|a_{ij}| \le M$, $a_{ij}$ – выигрыш игрока А (проигрыш игрока Б) в случае когда игрок А выбрал стратегию $i$, а игрок Б стратегию $j$. Отождествим себя с игроком Б. И предположим, что игра повторяется $N \gg 1$ раз (это число может быть заранее неизвестно). Мы находимся в условиях примера 4 с $\lambda(\omega^k, \zeta_j^k) = \sum_{i=1}^{n} \omega_i^k a_{ij}$, то есть

$$f_k(x) = \langle \omega^k, Ax \rangle, \ x \in S_n(1),$$

где $\omega^k$ – вектор (вообще говоря, зависящий от всей истории игры до текущего момента включительно, в частности, как-то зависящий и от текущей стратегии (не хода) игрока Б, заданной распределением вероятностей (результат текущего разыгрывания (ход Б) игроку А не известен)) со всеми компонентами равными 0, кроме одной компоненты, соответствующей ходу А на шаге $k$, равной 1. Хотя функция $f_k(x)$ определена на единичном симплексе, по "правилам игры" вектор $x^k$ имеет ровно одну единичную компоненту, соответствующую ходу Б на шаге $k$, остальные компоненты равны нулю. Обозначим цену игры

$$C = \max_{\omega \in S_n(1)} \min_{x \in S_n(1)} \langle \omega, Ax \rangle = \min_{x \in S_n(1)} \max_{\omega \in S_n(1)} \langle \omega, Ax \rangle. \text{ (теорема фон Неймана о минимаксе)}$$

**Замечание 11.** Отметим, что с помощью онлайн оптимизации и экспоненциального взвешивания можно похожим образом проинтерпретировать и вариант теоремы о минимаксе для векторнозначной функции выигрыша – теорему Блэкуэлла о достижимости [9], [27], которая используется, например, при построении калибруемых предсказаний.

Пару векторов $(\omega, x)$, доставляющих решение этой минимаксной задачи (т.е. седловую точку), назовем равновесием Нэша. По определению (это неравенство восходит к Ханнану [9], [27])

$$\min_{x \in S_n(1)} \frac{1}{N} \sum_{k=1}^{N} f_k(x) \le C.$$

Тогда, если мы (игрок Б) будем придерживаться рандомизированной стратегии МЗС2, выбирая $\{x^k\}$, то по теореме 2 с вероятностью $\geq 1-\sigma$ (в случае когда $N$ заранее известно оценку можно уточнить)

$$\frac{1}{N}\sum_{k=1}^{N}f_k(x^k) - \min_{x \in S_n(1)}\frac{1}{N}\sum_{k=1}^{N}f_k(x) \leq \frac{2M}{\sqrt{N}}\left(\sqrt{\ln n} + \sqrt{2\ln(\sigma^{-1})}\right),$$

т.е. с вероятностью $\geq 1-\sigma$ наши потери ограничены

$$\frac{1}{N}\sum_{k=1}^{N}f_k(x^k) \leq C + \frac{2M}{\sqrt{N}}\left(\sqrt{\ln n} + \sqrt{2\ln(\sigma^{-1})}\right).$$

Самый плохой для нас случай (с точки зрения такой оценки) – это когда игрок А тоже "знает" теорему 2, и действует (выбирая $\{\omega^k\}$) согласно МЗС2 (точнее версии МЗС2 для максимизации вогнутых функций на симплексе). Очевидно, что если и А и Б будут придерживаться МЗС2, то они сойдутся к равновесию Нэша (седловой точке), причем чрезвычайно быстро [26]:

$$\frac{8M\left(\ln n + 2\ln(\sigma^{-1})\right)}{\varepsilon^2} \text{ – итераций};$$

$$O\left(n + M\frac{s\ln n\left(\ln n + \ln(\sigma^{-1})\right)}{\varepsilon^2}\right) \text{ – общее число арифметических операций},$$

где $s \leq n$ – среднее число элементов в строках и столбцах матрицы $A$. Отсюда видно, что если $\varepsilon$ (зазор двойственности [3], [7]) – не очень малое, то может случиться, что общее число арифметических операций будет много меньше числа элементов матрицы $A$ отличных от нуля, в то время как любой детерминированный способ поиска равновесия Нэша потребовал бы прочтения как минимум половины элементов матрицы $A$ [6], [26]. Другой способ содержательной интерпретации описанного метода базируется на прямодвойственности МЗС [3].

**Замечание 12.** Если использовать методы работ [29], [30] то можно получить зависимость сложности от $\varepsilon$ вида $\varepsilon^{-1}$, но при этом число операций увеличится не менее чем в $n$ раз. Поскольку для таких задач вполне естественным является соотношение $\varepsilon \gg n^{-1}$, то выгоднее использовать описанный в этом примере алгоритм.

Есть основания полагать, что описанный здесь метод оптимален не только с точки зрения числа итераций, но и с точки зрения "стоимости" шага. Особенно ярко это проявляется, когда $s \ll n$ [26]. Нам кажется, что описанная в этом примере методология может оказаться полезной в huge-scale optimization (не только для рассмотренной здесь задачи). В частности, при изучении того, как оптимальным образом можно учитывать разреженность задачи [31] и как оптимально организовывать покомпонентный спуск [32]. Отметим в связи с вышесказанным другой пример задачи выпуклой оптимизации, когда удается получить ответ (с требуемой точностью), с большим запасом не просматривая весь объем имеющихся данных (общая проблема здесь: как понять, что просматривать, а что нет): метод наименьших квадратов с разреженной структурой [33]. Нам представляется, что именно эта ветвь более общего и бурно развивающегося в последнее время направления huge-scale оптимизации является сейчас одной из наиболее интересных, как с практической, так и с теоретической точки зрения, и основные результаты здесь еще впереди. Уточним, что речь идет о задачах, приходящих из: машинного обучения (в частности, распознавания изображений), моделирования различных сетей огромных размеров типа сети Интернет или транспортных сетей, биоинформатики, численных методов (проектирование конструкций методом конечных элементов) и ряда других приложений. Все эти задачи отличают колоссальные размеры. Скажем, для задачи ранжирования web-страниц необходимо решать вспомогательную оптимизационную задачу в пространстве, размерность которого больше миллиарда [26]. Но помимо размеров их отличает некоторые релаксированные требования к решению. Например, нет никакой необходимости ранжировать абсолютно все web-страницы по заданному запросу и делать это очень точно. Достаточно, чтобы качественно выдавались только первые сто наиболее значимых (больших) компонент ранжирующего вектора, причем важны не столько сами значения этих компонент, сколько их порядок. Также эти задачи довольно специальные, то есть использовать концепцию черного ящика [1] для оценки числа требуемых итераций, как правило, не представляется возможным. Более того, оценки имеются, в основном, только на число шагов (итераций). Поскольку общее время работы алгоритма определяется произведением числа итераций на стоимость одной итерации, то возникает игра между числом итераций и стоимостью итерации. Описанный в этом примере метод как раз играет в эту игру. А именно, он увеличивает за счет рандомизации число итераций в $\ln(\sigma^{-1})$ раз, при этом итерация становится в $n/\ln n$ раз дешевле. Наконец, важной составляющей многих задач является разреженная структура. Тут имеется два варианта: разреженная структура решения и данных. В первом случае (частично уже затронутым выше на примере ранжирования web-страниц) часто речь идет о подмене исходной задачи, вычислительно более привлекательной задачей, по решению которой можно получить приближенное решение исходной. Пожалуй, наиболее ярким примером здесь является сжатие измерений (см. [35], и цитированную там литературу). В случае разреженных данных представляется перспективным использование специальных покомпонентных спусков и исследование вычислительных особенностей пересчета различных классов функций многих переменных в случае изменения лишь не большого числа их аргументов [31], [32]. В заключение отметим, что важными составляющими анализа эффективности методов (в виду специфики описываемых задач) является исследование возможности распараллеливания [34] (в рассматриваемом нами примере 5 это возможно сделать [6], [26]) и вероятностный анализ в среднем [35] (или для почти всех входов). В отличие от Computer Science [36], [37], в численных методах выпуклой оптимизации такой анализ можно встретить не часто (все же кое-что есть, например, вероятностный анализ симплекс-метода Спилманом и Тэнгом [38]). Хотя уже сейчас (в связи с бурным развитием идей концентрации меры [37], [39]) становится все более и более ясно, что в пространствах огромных размеров такой анализ необходим, и может многое дать [35].

Отметим в заключение, что к рассмотренной в примере 5 антагонистической матричной игре (с $M=1$) сводится Google problem: задача поиска вектора Фробениуса–Перрона стохастической матрицы $P$ огромных размеров [26]: $A = P^T - I$. Разработанный в [26] и описанный в примере 5 подход позволяет учитывать разреженную структуру матрицы $P$, заменяя в алгоритме из [8] $n$, которое с точностью до логарифмических факторов входит линейно в оценку общей трудоемкости метода (то есть с учетом затрат на каждой итерации) на $s \ll n$.

## Литература